\documentclass[12pt]{amsart}

\newtheorem{theorem}{Theorem}[section]
\newtheorem{lemma}[theorem]{Lemma}

\theoremstyle{plain}
\newtheorem{definition}[theorem]{Definition}

\newtheorem{remark}[theorem]{Remark}

\newtheorem{question}[theorem]{Question}

\theoremstyle{definition}

\theoremstyle{remark}

\numberwithin{equation}{section}

\usepackage{fancyhdr}
\usepackage{amscd}
\usepackage{amsmath}
\usepackage{amsthm}
\usepackage{amssymb}

\begin{document}

\title{Tangent spaces to motivic cohomology groups \\ \vspace{3 mm}}

\author{Sen Yang}
\address{Yau Mathematical Sciences Center, 
Tsinghua University, 
Beijing,100084, China.}
\email{syang@math.tsinghua.edu.cn; senyangmath@gmail.com}




\maketitle

\begin{abstract}
By using Green-Griffiths' results on tangent spaces to algebraic cycles \cite{GGtangentspace}, we study the tangent space to $CH^{2}(X,1)$, where $X$ is a nonsingular projective curve over a field $k$ of characteristic $0$.
\end{abstract}

\tableofcontents

\section{\textbf{Introduction}}
\label{Introduction}
For a regular scheme $X$ over a field $k$, Voevodsky defines motivic cohomology $H^{p}_{M}(X, \mathbb{Z}(q))$ and proves the following identification(for $k$ a perfect field)
\begin{equation}
 H^{p}_{M}(X, \mathbb{Z}(q)) = CH^{q}(X,2q-p),
\end{equation}
where $CH^{q}(X,2q-p)$ is Bloch's higher Chow group.

However, this theory of motivic cohomology has one obvious deficiency: it does not take care of the infinitesimal structure which is very important for studying deformation problems. For example, for $X$ a regular scheme over a field $k$, the $j$-th infinitesimal thickening $X_{j}=(X, O_{X}[t]/(t^{j+1}))$, is  a typical non-reduced scheme. Motivic cohomology can't distinguish the difference between $X$ and $X_{j}$.

 As the first attempt to understand what motivic cohomology of infinitesimal thickenings might mean, Bloch-Esnault \cite{BlochEsnault} introduces an additive version of higher Chow groups. In \cite{BlochLectures}, Bloch poses the following question:
\begin{question}
Assuming one has a good definition of motivic cohomology(of infinitesimal thickenings), still denoted $H^{p}_{M}(X, \mathbb{Z}(q))$, what should  the \textbf{tangent space}  
\[
 TH^{p}_{M}(X, \mathbb{Z}(q)):= \mathrm{Ker} \{H^{p}_{M}(X \times \mathrm{Spec}(k[\varepsilon]), \mathbb{Z}(q)) \xrightarrow{\varepsilon=0} H^{p}_{M}(X, \mathbb{Z}(q)) \}
\]
mean? Here $\mathrm{Spec}(k[\varepsilon])$ denotes the dual number, $\varepsilon^{2}=0$.
\end{question}

According to identification (1.1), we rewrite this question in terms of higher Chow groups as follows:
\begin{question}
Can one give a good definition of higher Chow groups of infinitesimal thickenings? If so, what should  the \textbf{tangent space} 
\[
 TCH^{q}(X,p):= \mathrm{Ker} \{CH^{q}(X \times \mathrm{Spec}(k[\varepsilon]),p) \xrightarrow{\varepsilon=0}  CH^{q}(X,p)\}
\]
mean? 
\end{question}

In this note, we focus on \textbf{Question 1.2}. After a brief review on background in Section 2, we study the tangent space to $CH^{2}(X,1)$, where $X$ is a nonsingular projective curve, in Section 3.

\textbf{Notations and conventions}.
For any abelian group $M$, $M_{\mathbb{Q}}$ denotes the image of $M$ in $M \otimes_{\mathbb{Z}} \mathbb{Q}$. $\mathrm{Spec}(k[\varepsilon])$ denotes the dual number, $\varepsilon^{2}=0$. 
 
\section{\textbf{Background}}
\label{Background}
Let $X$ be a quasiprojective variety over a field $k$ of characteristic $0$. Higher Chow groups are defined by Bloch in \cite{BlochHigherChow}. Let $\bigtriangleup^{m}:= \mathrm{Spec}(k[t_{0}, . . . , t_{m}]/(\sum_{i=0}^{m} t_{i}-1))$ be the $m$-simplex, Bloch introduces 
\begin{align*}
Z^{p}(X,m)=\{& \xi \in Z^{p}(X \times \bigtriangleup^{m}) \mid \\
 & \xi  \ \mbox{meets all faces} \ \{t_{i_{1}}=\dots= t_{i_{l}}=0, l\geq 1 \} \ \mbox{properly} \}.
\end{align*}
Now set $\partial_{j}:Z^{p}(X,m) \to Z^{p}(X,m-1)$, the restriction to the $j$-th face, given by $t_{j}=0$. The boundary map 
\[
\partial=\sum_{j=0}^{m}(-1)^{j}\partial_{j}: Z^{p}(X,m) \to Z^{p}(X,m-1)
\]
satisfies $\partial^{2}=0$. So one obtains the following complex $(\mathcal{Z}^{p}(X,\bullet), \partial)$:
\[
\dots \to Z^{p}(X,m+1) \xrightarrow{\partial} Z^{p}(X,m) \xrightarrow{\partial} Z^{p}(X,m-1) \to \dots
\]

\begin{definition} \cite{BlochHigherChow}
The higher Chow group $CH^{p}(X,m)$ is defined to be the homology of  $(\mathcal{Z}^{p}(X,\bullet), \partial)$ at position $m$
\[
  CH^{p}(X,m):= H_{m}(\mathcal{Z}^{p}(X,\bullet), \partial).
\]
\end{definition}

Let $X$ be a nonsingular projective variety over a field $k$ of characteristic $0$ and let $K^{M}_{p}(O_{X})$ denote the Milnor K-theory sheaf associated to the presheaf
\[
  U \to  K^{M}_{p}(O_{X}(U))).
\]
According to \cite{MS}, one can relate the higher Chow groups $CH^{p}(X,m)$ with $H^{p-m}_{Zar}(X, K_{p}^{M}(O_{X}))$ for $0 \leq m \leq 2$ in the following diagram
\[
  \begin{CD}
   Z^{p}(X,2) @>\mathrm{N}>> \bigoplus\limits_{x \in X^{(p-2)}}K^{M}_{2}(k(x)) \\
   @VVV    @VV\mathrm{Tame}V\\
   Z^{p}(X,1) @>\mathrm{N}>> \bigoplus\limits_{x \in X^{(p-1)}}K^{M}_{1}(k(x))\\
   @VVV    @VV\mathrm{div}V\\
   Z^{p}(X) @>=>> \bigoplus\limits_{x \in X^{(p)}}K^{M}_{0}(k(x)).\\
  \end{CD}
\]

Here N stands for the Norm map, Tame and div are short for Tame symbol and divisor map respectively. This leads to the following identifications.
\begin{theorem}  \label{theorem:Bloch-type formulas}

 \cite{MS,MS2,Soule}  \textbf{Bloch-type formulas}.
Let $X$ be a nonsingular projective variety over a field $k$ of characteristic $0$.  One has the following identifications,  for  $0 \leq m \leq 2$,
\[
  CH^{p}(X,m)_{\mathbb{Q}}=H^{p-m}_{\mathrm{Zar}}(X, K^{M}_{p}(O_{X}))_{\mathbb{Q}}.
\]
\end{theorem}

 \begin{proof}
 When $m=0$, the above formula is Soul\'e's variant of Bloch-Quillen identification.  
 M\"{u}ller-Stach proves the cases of $m=1$ and $m=2$ in \cite{MS,MS2}. 
 \end{proof}

 Based on this theorem, the higher Chow groups $CH^{p}(X,m)$, $0\leq m \leq 2$, can be described as follows, see \cite{Lewislecture}. 

1.) $CH^{p}(X,0)=CH^{p}(X)$, where $CH^{p}(X)$ is the classical Chow group of algebraic cycles modulo rational equivalence.

2.)  $CH^{p}(X,1)$ is represented as a quotient:
 \[
  CH^{p}(X,1)= \dfrac{\mathrm{Ker(div)}}{\mathrm{Im(Tame)}}.
 \]
In other words, an element of $CH^{p}(X,1)$ is of the form 
\[
\{\sum_{j}\{Y_{j}, f_{j}\mid_{Y_{j}}\}\mid
\mbox{codim}(Y_{j})=p-1, f_{j} \in k(Y_{j})^{*} \mbox{and} \sum_{j} \mathrm{div}(f_{j})=0 \}
\]
and modulo the image of Tame symbol.

3.) $CH^{p}(X,2)$ is represented by classes in the kernel of Tame symbol, modulo 
the image of higher Tame symbol.

Let $TK^{M}_{p}(O_{X})$ denote the tangent space to $K^{M}_{p}(O_{X})$, it is known that $TK^{M}_{p}(O_{X}) \cong \Omega_{X/ \mathbb{Q}}^{p-1}$. By imitating the infinitesimal method of Bloch, one can define tangent space via the above Bloch-type formulas:
\begin{definition}  \cite{BlochLectures} \label{definition: tangentspaceagree}
Let $X$ be a nonsingular projective variety over a field $k$ of characteristic $0$.  For  $0 \leq m \leq 2$, the tangent space to $CH^{p}(X,m)$ is defined to be 
\[
 TCH^{p}(X,m) := H^{p-m}(X, TK^{M}_{p}(O_{X})) = H^{p-m}(X, \Omega_{X/ \mathbb{Q}}^{p-1}).
\]
\end{definition}

One natural question is: does this formal definition carry any concrete geometric meaning?

\section{\textbf{Tangent spaces to $CH^{2}(X,1)$}}
\label{Tangent spaces to $CH^{2}(X,1)$}
In this section, to fix notations, $X$ is a nonsingular projective curve over a field $k$ of characteristic $0$.  
We shall give a concrete geometric meaning to the tangent space $TCH^{2}(X,1)$ by using Green-Griffiths' results in \cite{GGtangentspace}. The key to our approach is the following splitting commutative diagram.

\begin{theorem} \cite{Y-2}  \label{theorem:maintheorem}
For $q=2$ and $j=1$ in $\mathrm{Theorem \ 3.12}$ of \cite{Y-2}, we have the following splitting commutative diagram in which the Zariski sheafification of each column is a flasque resolution of \ $\Omega_{X/ \mathbb{Q}}^{1}$,  $K_{2}(O_{X}[\varepsilon]) $ and $ K_{2}(O_{X})$ respectively:
\[\displaystyle
  \begin{CD}
     \Omega_{k(X)/ \mathbb{Q}}^{1} @<\mathrm{Chern-1}<<  K_{2}(k(X)[\varepsilon]) @>\varepsilon=0>> K_{2}(k(X)) \\
     @V \partial^{0,-2}_{1}VV @Vd^{0,-2}_{1, X_{1}}VV @Vd^{0,-2}_{1,X}V\mathrm{Tame}V\\
     \bigoplus\limits_{x \in X^{(1)}} H_{x}^{1}(\Omega_{X/\mathbb{Q}}^{1}) @<\mathrm{Chern-2}<<  \bigoplus\limits_{x[\varepsilon] \in X[\varepsilon]^{(1)}}K_{1}(O_{X,x}[\varepsilon] \ \mathrm{on} \ x[\varepsilon]) @>\varepsilon=0>> \bigoplus\limits_{x \in X^{(1)}} K_{1}(k(x)) \\
     @VVV @VVV @VVV\\
      0 @. 0 @. 0,
  \end{CD}
\]
where $\mathrm{Chern-1}$ and $\mathrm{Chern-2}$ are induced by Chern character from K-theory to negative cyclic homology.
\end{theorem}

Since $K_{2}(k(X)[\varepsilon])$ can be identified with Milnor K-group $K^{M}_{2}(k(X)[\varepsilon])$(and similarly for $K_{2}(k(X))$ and $K_{1}(k(x))$,  the above commutative diagram can be rewritten as:
\[\displaystyle
  \begin{CD}
     \Omega_{k(X)/ \mathbb{Q}}^{1} @<\mathrm{Chern-1}<<  K^{M}_{2}(k(X)[\varepsilon]) @>\varepsilon=0>> K^{M}_{2}(k(X)) \\
     @V \partial^{0,-2}_{1}VV @Vd^{0,-2}_{1, X_{1}}VV @Vd^{0,-2}_{1,X}V\mathrm{Tame}V\\
     \bigoplus\limits_{x \in X^{(1)}} H_{x}^{1}(\Omega_{X/\mathbb{Q}}^{1}) @<\mathrm{Chern-2}<<  \bigoplus\limits_{x[\varepsilon] \in X[\varepsilon]^{(1)}}K_{1}(O_{X,x}[\varepsilon] \ \mathrm{on} \ x[\varepsilon]) @>\varepsilon=0>> \bigoplus\limits_{x \in X^{(1)}} K^{M}_{1}(k(x)) \\
     @VVV @VVV @VVV\\
      0 @. 0 @. 0.
  \end{CD}
\]

In the following,  we shall describe $K_{1}(O_{X,x}[\varepsilon] \ \mathrm{on} \ x[\varepsilon])$,  Chern-1 and Chern-2 explicitly.

\subsection{Green-Griffiths' $Arcs$  and Nenashev's result}
\label{Green-Griffiths' Arcs  and Nenashev's result}

Inspired by Gillet and Grayson's work on attaching a simplicial set to any exact category, Nenashev in \cite{Nenashev} provides a way to describe $K_{1}$ of any exact category in terms of generators and relations. 
Let $\mathcal{E}$ denote an exact category in the following.
\begin{definition} \cite{Nenashev}
 A double short exact sequence in $\mathcal{E}$ is a pair of short exact sequences on the same objects:
\begin{equation}
\begin{cases}
 \begin{CD}
   0 \to A \xrightarrow{f_{1}} B \xrightarrow{g_{1}}   C \to 0\\
   0 \to A \xrightarrow{f_{2}} B \xrightarrow{g_{2}}   C \to 0.
 \end{CD}
\end{cases}
\end{equation} 

\end{definition}
In particular, if $A \in \mathcal{E}$ and $\alpha \in \mathrm{Aut}(A)$, we can associate a double short exact sequence to $\alpha$:
\begin{equation}
\begin{cases}
 \begin{CD}
   0 \to A \xrightarrow{1} A \to 0\\
   0 \to A \xrightarrow{\alpha} A \to 0.
 \end{CD}
\end{cases}
\end{equation} 

Now one defines an abelian group generated by these double short exact sequences.
\begin{definition} \cite{Nenashev} \label{definition: TwoRelation}
 We define $A(\mathcal{E})$ to be the abelian group generated by all double short exact sequences subjecting to the following two relations:
 \begin{itemize}
  \item (1). If $f_{1}=f_{2}$ and  $g_{1}=g_{2}$, then the double short exact sequence is $0$.
  \item (2). $3\times3$ relations defined in $\mathrm{Proposition \ 2.1}$ in \cite{Nenashev}.
 \end{itemize}
\end{definition}

The main theorem in \cite{Nenashev} by Nenashev says:
\begin{theorem}\cite{Nenashev} \label{theorem: Nenashev}
 For any exact category $\mathcal{E}$, there is an isomorphism between the following two abelian groups:
\[
 A(\mathcal{E}) \xrightarrow{\simeq} K_{1}(\mathcal{E}).
\]
\end{theorem}

The following definition is used by Green-Griffiths \cite{GGtangentspace}, see chap 6, page 68:
\begin{definition} \cite{GGtangentspace}
 Let $X$ be a nonsingular projective curve over a field $k$ of characteristic $0$. For a closed point $x \in X$ with local uniformizer $f$, let $\mathrm{div}(f+\varepsilon f_{1}) $ denote $\mathrm{Spec}(O_{X,x}[\varepsilon]/(f+\varepsilon f_{1}))$, where $f_{1} \in O_{X,x}$.  An arc is defined to be a pair of the form 
\[
 \{\mathrm{div}(f+\varepsilon f_{1}), g+ \varepsilon g_{1}\mid_{\mathrm{div}(f+\varepsilon f_{1})} \},
\]
where $g \in k(x)^{\ast}$ and $g+ \varepsilon g_{1} \in (k(x)[\varepsilon])^{\ast}$ and furthermore we assume $\mathrm{div}(f+\varepsilon f_{1})\cap \mathrm{div}(g+ \varepsilon g_{1}) = \emptyset$.
\end{definition}

In the following, we use $Arcs$ to denote the set of $arc$: $Arcs= \{ arc \}$.
\begin{remark}
Considering $\mathrm{Spec}(k[\varepsilon])=\mathrm{Spec}(k[t]/(t^2))$, one can write the above $arc$ as the form
\[
 \{\mathrm{div}(f+tf_{1}), g+tg_{1}\mid_{\mathrm{div}(f+tf_{1})} \}.
\]
This is the notation used in \cite{GGtangentspace}. Intuitively, one can think of $\mathrm{div}(f+tf_{1})$(resp. $g+tg_{1}$) as the $1^{st}$ order deformation of $\mathrm{div}(f)$(resp. $g$). 
\end{remark}

We want to identify the above $arc$ as an element of $K_{1}(O_{X,x}[\varepsilon] \ \mathrm{on} \ x[\varepsilon])$, so we need the following theorem, $Exercise \ 5.7$ of Thomason-Trobaugh \cite{TT}. 

\begin{theorem} \cite{TT}
Let $X$ be a scheme with an ample family
of line bundles. Let $i : Y \to X$ be a regular closed immersion ([SGA 6]
VII Section 1) defined by ideal $J$. Suppose $Y$ has codimension $k$ in X.
Then $K(X \ \mathrm{on} \ Y)$ is homotopy equivalent to the Quillen K-theory
of the exact category of pseudo-coherent $O_{X}$-modules supported
on the subspace $Y$ and of Tor-dimension $\leq k$ on $X$.
\end{theorem}

According to this theorem, $K_{1}(O_{X,x}[\varepsilon] \ \mathrm{on} \ x[\varepsilon])$ can be considered as a K-group of the exact category of pseudo-coherent $O_{X,x}[\varepsilon]$-modules supported
on the subspace $x[\varepsilon]$ and of Tor-dimension $\leq 1$ on $O_{X,x}[\varepsilon]$.
$O_{X,x}[\varepsilon]/(f+\varepsilon f_{1}))$ is such a module.

Considering $g+ \varepsilon g_{1}$ as an automorphism of $O_{X,x}[\varepsilon]/(f+\varepsilon f_{1})$, we can associate a double short exact sequence  to $\{\mathrm{div}(f+\varepsilon f_{1}), g+ \varepsilon g_{1}\mid_{\mathrm{div}(f+\varepsilon f_{1})} \}$ 
\begin{equation}
\begin{cases}
 \begin{CD}
   0 \to O_{X,x}[\varepsilon]/(f+\varepsilon f_{1}) \xrightarrow{1} O_{X,x})[\varepsilon]/(f+\varepsilon f_{1}) \to 0 \\
   0 \to O_{X,x}[\varepsilon]/(f+\varepsilon f_{1}) \xrightarrow{g+ \varepsilon g_{1}} O_{X,x}[\varepsilon]/(f+\varepsilon f_{1}) \to 0.
 \end{CD}
\end{cases}
\end{equation} 
According to Theorem ~\ref{theorem: Nenashev}, this double short exact sequence is an element of $K_{1}(O_{X,x}[\varepsilon] \ \mathrm{on} \ x[\varepsilon])$.  The above discussion shows:
\begin{lemma} \label{lemma: arcK1}
An $\{\mathrm{div}(f+\varepsilon f_{1}), g+ \varepsilon g_{1}\mid_{\mathrm{div}(f+\varepsilon f_{1})} \}$ defines an element of $K_{1}(O_{X,x}[\varepsilon] \ \mathrm{on} \ x[\varepsilon])$. 
\end{lemma}

\subsection{Green-Griffiths' tangent maps and Chern character}
\label{Green-Griffiths' tangent maps and Chern character}
In this subsection, we recall Green-Griffiths' geometric descriptions of tangent maps from $Arcs$ to local cohomology groups, and compare Green-Griffiths' tangent maps with Chern character maps in Theorem ~\ref{theorem:maintheorem}. 

\textbf{Describing Chern-1}. An element of $K^{M}_{2}(k(X)[\varepsilon])$ is given by Steinberg symbol $\{f+ \varepsilon f_{1}, g+ \varepsilon g_{1} \}$, where $f,g \in k(X)^{*}$ and $f_{1},g_{1} \in k(X)$. 
According to \cite{Lodaycyclic}-Section 8.4(page 275), the following composition, call it Ch,
\[
 \mathrm{Ch}: K^{M}_{2}(k(X)[\varepsilon]) \to K_{2}(k(X)[\varepsilon]) \xrightarrow{\mathrm{Chern}} HN_{2}(k(X)[\varepsilon])  \to \Omega^{2}_{k(X)[\varepsilon] / \mathbb{Q}}
\]
sends Steinberg symbol $\{s, t \}$ to $s^{-1}t^{-1}dsdt$. Furthermore, there is a truncation map from $\Omega^{2}_{k(X)[\varepsilon] / \mathbb{Q}}$ to $\Omega^{1}_{k(X)/ \mathbb{Q}}$
\begin{align*}
\dfrac{\partial}{\partial \varepsilon}\mid_{\varepsilon=0}: & \hspace{8mm} \Omega^{2}_{k(X)[\varepsilon] / \mathbb{Q}} \longrightarrow \Omega^{1}_{k(X)/ \mathbb{Q}}\\
& (a+b\varepsilon)d(x+y\varepsilon)d(z+w\varepsilon) \to a(wdx-ydz).
\end{align*}
 
It is classical that the Chern-1 map in Theorem ~\ref{theorem:maintheorem},  which can be described as the composition of Ch and the truncation $\dfrac{\partial}{\partial \varepsilon}\mid_{\varepsilon=0}$,  is of the following form:
\begin{align*}
\mathrm{Chern}-1: \ & K^{M}_{2}(k(X)[\varepsilon]) \xrightarrow{\mathrm{Ch}} \Omega^{2}_{k(X)[\varepsilon] / \mathbb{Q}} \xrightarrow{\dfrac{\partial}{\partial \varepsilon}\mid_{\varepsilon=0}} \Omega^{1}_{k(X)/ \mathbb{Q}}\\
& \{f+ \varepsilon f_{1}, g+ \varepsilon g_{1} \} \longrightarrow \dfrac{g_{1}df-f_{1}dg}{fg}.
\end{align*}

This is the form used by Green-Griffiths \cite{GGtangentspace}, see page 130.

\textbf{Green-Griffiths' tangent map}.  Next,  we describe the Chern-2 map in Theorem ~\ref{theorem:maintheorem}.
We begin with recalling Green-Griffiths' description of $\partial^{0,-2}_{1}$, see \cite{GGtangentspace}-page 105 for more discussions(for $\mathbb{C}(X) \to \bigoplus\limits_{y \in X^{(1)}} H_{y}^{1}(O_{X})$ and works similarly here). 
Working locally in a Zariski open affine neighborhood U, we can write an element $ \beta \in \Omega_{K(X) / \mathbb{Q}}^{1}$ as 
\[
 \beta = \dfrac{h \ dg}{f^{l_{1}}_{1}\dots f^{l_{k}}_{k}},
\]
where $f_{1}, \dots, f_{k}, h, g \in \Gamma(U,O_{U})$ are relatively prime and $f_{i}'s$ are irreducible. Set $x_{i}=\{f_{i} = 0 \}$ and let $\beta_{i}$ denote 
\[
\beta_{i} = \dfrac{h \ dg}{f^{l_{1}}_{1}\dots \hat{f}^{l_{i}}_{i} \dots f^{l_{k}}_{k}},
\]
where $\hat{f}^{l_{i}}_{i}$ means to omit the $i^{th}$ term. By abuse of notations, we still use $\beta_{i}$ to denote the following diagram, where $R= \Gamma(U,O_{U})$,
\begin{equation}
\begin{cases}
 \begin{CD}
   R_{(f^{l_{i}}_{i})} @>f^{l_{i}}_{i}>> R_{(f^{l_{i}}_{i})} @>>> R_{(f^{l_{i}}_{i})}/(f^{l_{i}}_{i}) @>>> 0  \\
   R_{(f^{l_{i}}_{i})} @>\beta_{i}>> \Omega_{R_{(f^{l_{i}}_{i})}/ \mathbb{Q}}^{1}.
 \end{CD}
\end{cases}
\end{equation}

\begin{lemma} \cite{GGtangentspace}-$\mathrm{page \ 105}$. With above notations, $\partial^{0,-2}_{1}$ can be described as follows:
\begin{align*}
\partial^{0,-2}_{1}: & \  \Omega_{K(X)/ \mathbb{Q}}^{1} \to \bigoplus_{x \in X^{(1)}} H_{x}^{1}(\Omega_{X/\mathbb{Q}}^{1})\\
& \beta \longrightarrow \sum_{i}\beta_{i}.
\end{align*}
\end{lemma}

Now we recall Green-Griffiths' description of $d^{0,-2}_{1, X_{1}}$. An element of $K^{M}_{2}(k(X)[\varepsilon])$ is given by Steinberg symbol $\{f+ \varepsilon f_{1}, g+ \varepsilon g_{1} \}$, where $f,g \in k(X)^{*}$ and $f_{1},g_{1} \in k(X)$. Working locally in a Zariski open affine neighborhood U, we write $\{f+ \varepsilon f_{1}, g+ \varepsilon g_{1} \}$ as a product of symbols of the form $\{a+ \varepsilon a_{1}, b+ \varepsilon b_{1} \}$ or its inverse, where $a,b,a_{1},b_{1} \in R=\Gamma(U,O_{U})$, $a \neq 0$ and $b \neq 0$.

For simplicity, we assume an element of $K^{M}_{2}(k(X)[\varepsilon])$ is given by Steinberg symbol $\{f+ \varepsilon f_{1}, g+ \varepsilon g_{1} \}$, where $f,g,f_{1},g_{1} \in \Gamma(U,O_{U})$, $f \neq 0$ and $g \neq 0$. 

\begin{lemma} \cite{GGtangentspace}-$\mathrm{page \ 130}$.
The differential 
\[
d^{0,-2}_{1,X_{1}}:  K^{M}_{2}(k(X)[\varepsilon]) \to \bigoplus_{x[\varepsilon] \in X[\varepsilon]^{(1)}}K_{1}(O_{X,x}[\varepsilon] \ \mathrm{on} \ x[\varepsilon])
\]
can be described as follows, 
\begin{align*}
  \{f+ \varepsilon f_{1}, g+ \varepsilon g_{1} \} \to & \{\mathrm{div}(f+ \varepsilon f_{1}), g+ \varepsilon g_{1} \mid_{\mathrm{div}(f+ \varepsilon f_{1})} \} \\
 & - \{\mathrm{div}(g+ \varepsilon g_{1}), f+ \varepsilon f_{1} \mid_{\mathrm{div}(g+ \varepsilon g_{1})}\}, 
\end{align*}
where $\mathrm{div}(f+ \varepsilon f_{1})$ denotes $\mathrm{Spec}(R_{(f)}[\varepsilon]/(f+ \varepsilon f_{1}))$ and we assume $\mathrm{div}(f+ \varepsilon f_{1}) \cap \mathrm{div}(g+ \varepsilon g_{1}) = \emptyset$. Otherwise, the image is defined to be $0$. Here, $\{\mathrm{div}(f+ \varepsilon f_{1}), g+ \varepsilon g_{1} \mid_{\mathrm{div}(f+ \varepsilon f_{1})} \} $ is considered as an element of $\bigoplus\limits_{x[\varepsilon] \in X[\varepsilon]^{(1)}}K_{1}(O_{X,x}[\varepsilon] \ \mathrm{on} \ x[\varepsilon])$ via $\mathrm{Lemma \ ~\ref{lemma: arcK1}}$.
\end{lemma}

Green-Griffiths \cite{GGtangentspace} defines a map, call it tangent , from the \textbf{set} of $Arcs$ to local cohomology groups(for $X$ a surface which can be easily  adopted to the case of $X$ a curve):
\[
 \mathrm{ tangent}:   Arcs   \ \longrightarrow \bigoplus_{x \in X^{(1)}}H_{x}^{1}(\Omega_{X/\mathbb{Q}}^{1}). 
\] 
\begin{definition} \cite{GGtangentspace}-$\mathrm{page \ 127}$.
 Recall that an element of $Arcs$ is defined to be pairs of the form 
\[
 \{\mathrm{div}(f+\varepsilon f_{1}), g+ \varepsilon g_{1}\mid_{\mathrm{div}(f+\varepsilon f_{1})} \},
\]
where $f_{1} \in O_{X,x}$, $g+ \varepsilon g_{1} \in (k(x)[\varepsilon])^{\ast}$ and furthermore we assume $\mathrm{div}(f+\varepsilon f_{1})\cap \mathrm{div}(g+ \varepsilon g_{1}) = \emptyset$.  

For the element $\{\mathrm{div}(f+\varepsilon f_{1}), g+\varepsilon g_{1} \mid_{\mathrm{div}(f+\varepsilon f_{1})} \}$, the following diagram
\begin{equation}
\begin{cases}
 \begin{CD}
   O_{X,x} @>f>> O_{X,x} @>>> O_{X,x}/(f) @>>> 0  \\
   O_{X,x} @>\frac{g_{1} df}{g}-\frac{f_{1} dg}{g}>> \Omega_{O_{X,x}/ \mathbb{Q}}^{1}
 \end{CD}
\end{cases}
\end{equation} 
gives an element $\alpha$ in $Ext_{O_{X,x}}^{1}(O_{X,x}/(f),\Omega_{O_{X,x}/ \mathbb{Q}}^{1})$. Noting that 
\[
H_{x}^{1}(\Omega_{X/\mathbb{Q}}^{1})=\varinjlim_{n \to \infty}Ext_{O_{X,x}}^{1}(O_{X,x}/(f)^{n},\Omega_{O_{X,x}/ \mathbb{Q}}^{1}),
\]
the image $[\alpha]$ of $\alpha$ under the limit is in $H_{x}^{1}(\Omega_{X/\mathbb{Q}}^{1})$ and it is defined to be the image of $\{\mathrm{div}(f+\varepsilon f_{1}), g+\varepsilon g_{1} \mid_{\mathrm{div}(f+\varepsilon f_{1})} \}$ under the $\mathrm{tangent}$ map.
\end{definition}

By abuse of notations, we still use $Arcs $ to denote the \textbf{subgroup} generated by the elements of $Arcs$ in $\bigoplus\limits_{x[\varepsilon] \in X[\varepsilon]^{(1)}} K_{1}(O_{X,x}[\varepsilon] \ \mathrm{on} \ x[\varepsilon])$. One can check that the above Green-Griffiths' tangent map factors out the two relations in Definition ~\ref{definition: TwoRelation}, so it is well-defined on the subgroup $Arcs$:
\[
 \mathrm{tangent }:   Arcs  \ \longrightarrow \bigoplus_{x \in X^{(1)}}H_{x}^{1}(\Omega_{X/\mathbb{Q}}^{1}). 
\]

One can easily check that the following diagram is commutative:
\[
  \begin{CD}
     \Omega_{k(X)/ \mathbb{Q}}^{1} @<\mathrm{Chern-1}<<  K^{M}_{2}(k(X)[\varepsilon]) @>\varepsilon=0>> K^{M}_{2}(k(X)) \\
     @V \partial^{0,-2}_{1}VV @Vd^{0,-2}_{1, X_{1}}VV @Vd^{0,-2}_{1,X}VV\\
     \bigoplus\limits_{x \in X^{(1)}} H_{x}^{1}(\Omega_{X/\mathbb{Q}}^{1}) @<\mathrm{tangent}<<   Arcs @>\varepsilon=0>> \bigoplus\limits_{x\in X^{(1)}} K^{M}_{1}(k(x)). \\ 
  \end{CD}
\]

Recall that we also have the following commutative diagram in Theorem ~\ref{theorem:maintheorem}
\[
  \begin{CD}
     \Omega_{k(X)/ \mathbb{Q}}^{1} @<\mathrm{Chern-1}<<  K^{M}_{2}(k(X)[\varepsilon]) @>\varepsilon=0>> K^{M}_{2}(k(X)) \\
     @V \partial^{0,-2}_{1}VV @Vd^{0,-2}_{1,X_{1}}VV @Vd^{0,-2}_{1,X}VV\\
     \bigoplus\limits_{x \in X^{(1)}} H_{x}^{1}(\Omega_{X/\mathbb{Q}}^{1}) @<\mathrm{Chern-2}<<  \bigoplus\limits_{x[\varepsilon] \in X[\varepsilon]^{(1)}} K_{1}(O_{X,x}[\varepsilon] \ \mathrm{on} \ x[\varepsilon]) @>\varepsilon=0>> \bigoplus\limits_{x\in X^{(1)}} K^{M}_{1}(k(x)). \\ 
  \end{CD}
\]

We have the following theorem:
\begin{theorem}
The $\mathrm{Chern}-2$ map agrees with Green-Griffiths' $\mathrm{tangent}$ map(on the image of $d^{0,-2}_{1,X_{1}}$).
\end{theorem}

Guided by the above infinitesimal study, we propose the higher chow group of the infinitesimal thickening as follows:
\begin{definition}
Let $X$ be a nonsingular projective curve over a field $k$ of characteristic $0$, 
the higher chow group of the infinitesimal thickening $CH^{2}(X[\varepsilon],1)$ is defined to 
\[
CH^{2}(X[\varepsilon],1):= \dfrac{\bigoplus\limits_{x[\varepsilon] \in X[\varepsilon]^{(1)}} K_{1}(O_{X,x}[\varepsilon] \ \mathrm{on} \ x[\varepsilon])}{\mathrm{Image\ of } \ d^{0,-2}_{1,X_{1}}}.
\]
\end{definition}

We conclude that the above approach is an experimental attempt, a full answer to Bloch's \textbf{Question 1.2} still remains to search.

\textbf{Acknowledgements}
The author must record that the main results in this note are due to Mark Green and Phillip Griffiths \cite{GGtangentspace}. This note is to interpret their ideas in a slight different language.

The author is very grateful to: Christophe Soul\'e for precious help on K-theory and for telling \cite{BlochEsnault} which started this note; to Ben Dribus, Jerome W.~Hoffman, James D. ~Lewis and Chao Zhang for discussions; to Stefan M\"{u}ller-Stach, Alexander Nenashev and anonymous comments on previous versions; to Joseph Steenbrink who helped polish the previous versions which greatly improved this note.

\end{document}